\documentclass[english,11pt]{article}
\usepackage{amsthm}
\usepackage[latin9]{inputenc}
\usepackage{babel}
\usepackage[hmargin=0.9in]{geometry}

\usepackage{subfig}
\usepackage{graphicx}
\usepackage[colorlinks=true,citecolor=blue,urlcolor=blue]{hyperref}

\usepackage{float}
\usepackage{url}
\usepackage{bm}

\usepackage{amsmath}
\usepackage{amssymb,amsfonts}
\usepackage{upgreek}
\usepackage[textsize=footnotesize]{todonotes}
\usepackage{stmaryrd}



\graphicspath{{./figures/}}

\newcommand{\order}{{\mathcal O}}
\newcommand{\yhat}{\hat{y}}
\newcommand{\Khat}{\hat{K}}
\newcommand{\rhohat}{\hat{\rho}}
\def\avg[#1]{\left\langle #1 \right\rangle}
\def\intavg[#1]{\llbracket #1 \rrbracket}
\def\ceff{c_\text{eff}}
\def\bk{\bf k}

\begin{document}
\title{Two-dimensional wave propagation in layered periodic media}
\author{Manuel Quezada de Luna\thanks{Department of Mathematics, Texas A\&M University. College Station, Texas
77843, USA. E-mail: \url{mquezada@math.tamu.edu}} \and David I. Ketcheson\thanks{King Abdullah University of Science and Technology (KAUST), CEMSE Division. E-mail: \url{david.ketcheson@kaust.edu.sa}}}
\maketitle

\begin{abstract}
We study two-dimensional wave propagation in materials whose properties
vary periodically in one direction only.  
High order homogenization is carried out to derive a dispersive effective
medium approximation.  One-dimensional materials with constant
impedance exhibit no effective dispersion.
We show that a new kind of effective dispersion may arise in
two dimensions, even in materials with constant impedance.
This dispersion is a macroscopic effect of microscopic diffraction caused by
spatial variation in the sound speed.
We analyze this dispersive effect by using high-order homogenization to derive
an anisotropic, dispersive effective medium.  We generalize to two dimensions
a homogenization approach that has been used previously for one-dimensional problems.
Pseudospectral solutions of the effective medium equations agree
to high accuracy with finite volume direct numerical simulations of the
variable-coefficient equations.
\end{abstract}

\section{Introduction}
Consider the propagation of acoustic waves in a two-dimensional medium whose
properties vary in one coordinate direction (say, $y$).  Such waves are described
by the PDE
\begin{align} \label{wave-equation}
    p_{tt} & = K(y) \nabla \cdot \left(\frac{1}{\rho(y)} \nabla p\right).
\end{align}
Here $p = p(x,y,t)$ is the pressure, $K(y)$ is the bulk modulus, and $\rho(y)$
is the material density.  We focus on the initial value problem in an unbounded spatial
domain.  We are interested in materials whose spatial variation is periodic:
\begin{align*}
K(y+\Omega) & = K(y), & \rho(y+\Omega) = \rho(y).
\end{align*}
Here $\Omega$ denotes the period.  In all numerical experiments and plots,
we set $\Omega=1$.

\begin{figure}
\begin{center}
\includegraphics[width=0.5\textwidth]{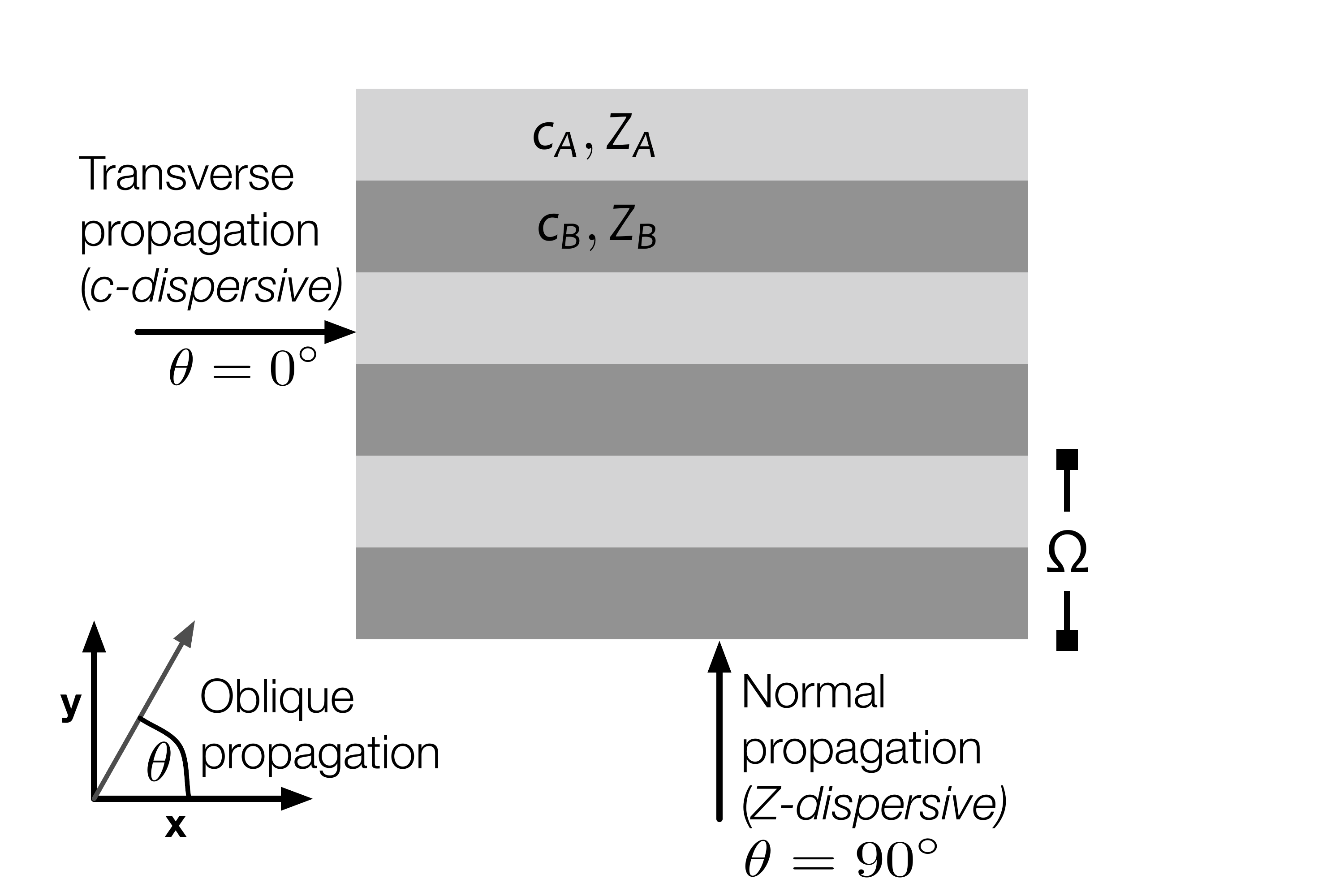}
\caption{Wave propagation in a layered periodic medium.
The material shown is repeated periodically in $y$ and extends 
infinitely in both coordinate directions.
A piecewise-constant medium is shown for simplicity, but arbitrary
periodic variation (in $y$) is considered.
The terms {\em normal} and {\em transverse} are used throughout this
work to denote propagation in the $y$ and $x$ coordinate directions, respectively.\label{fig:layered}}
\end{center}
\end{figure}

A simple example of such a medium is shown in Figure \ref{fig:layered}.
We refer to these as {\em layered} materials, though the coefficients
need not be piecewise-constant.  In subsequent sections, we frequently
use the terms {\em normal propagation} and {\em transverse propagation}
to refer to propagation normal to or parallel to the axis of homogeneity,
respectively (see Figure \ref{fig:layered}).

We consider the propagation of waves with characteristic wavelength $\lambda$
over a distance $L$ in a periodic medium with period $\Omega$ where
\[
\Omega < \lambda \ll L.
\]
Because the wavelength $\lambda$ is larger than the material period $\Omega$, the
waves ``see'' the medium as nearly homogeneous and travel at an effective velocity
related to averages of the material properties.  The study of wave propagation in
this regime has been the subject of much study; see
\cite{conca1997homogenization,fouque2007} and references therein.  
Many works focus exclusively on the lowest-order terms in the homogenized
equations.  The potential for dispersive higher-order terms
due to material periodicity was derived using Bloch expansions in \cite{santosa1991},
and computed explicitly for the case of a one-dimensional layered medium.
Later works have also studied effective dispersion in one-dimensional
periodic media, and further developed the relevant high-order homogenization techniques
for time-dependent problems \cite{fogarty1999,chen2001dispersive,Fish2001,Yong}.

The motivation for the present work comes from the discovery of a new kind of
effective dispersion, which is described and demonstrated briefly in the rest
of the introduction.
In Section \ref{sec:homog} we derive an effective medium approximation
to the variable-coefficient wave equation \eqref{wave-equation}.
The technique we use
is an extension of those appearing in \cite{chen2001dispersive,Yong,leveque2003}.
This seems to be the first explicit application of such high-order homogenization
to multidimensional materials.
In Section \ref{sec:dispersion} we explore the dispersion relation implied by
the effective medium equations, showing that the effective medium is 
anisotropic and
dispersive.  In Section \ref{sec:waves}, we examine the complementary roles
played by variation in the impedance and the sound speed.  We also
validate our homogenized model by comparing pseudospectral solutions of the
effective medium equations with finite volume direct
simulations of the variable-coefficient wave equation \eqref{wave-equation}.

All code used for computations in this work, along with Mathematica worksheets
used to derive the homogenized equations, are available at
\url{http://github.com/ketch/effective_dispersion_RR}.

\subsection{Effective dispersion in a layered medium}
\begin{figure}
\begin{center}
\includegraphics[width=1.0\textwidth]{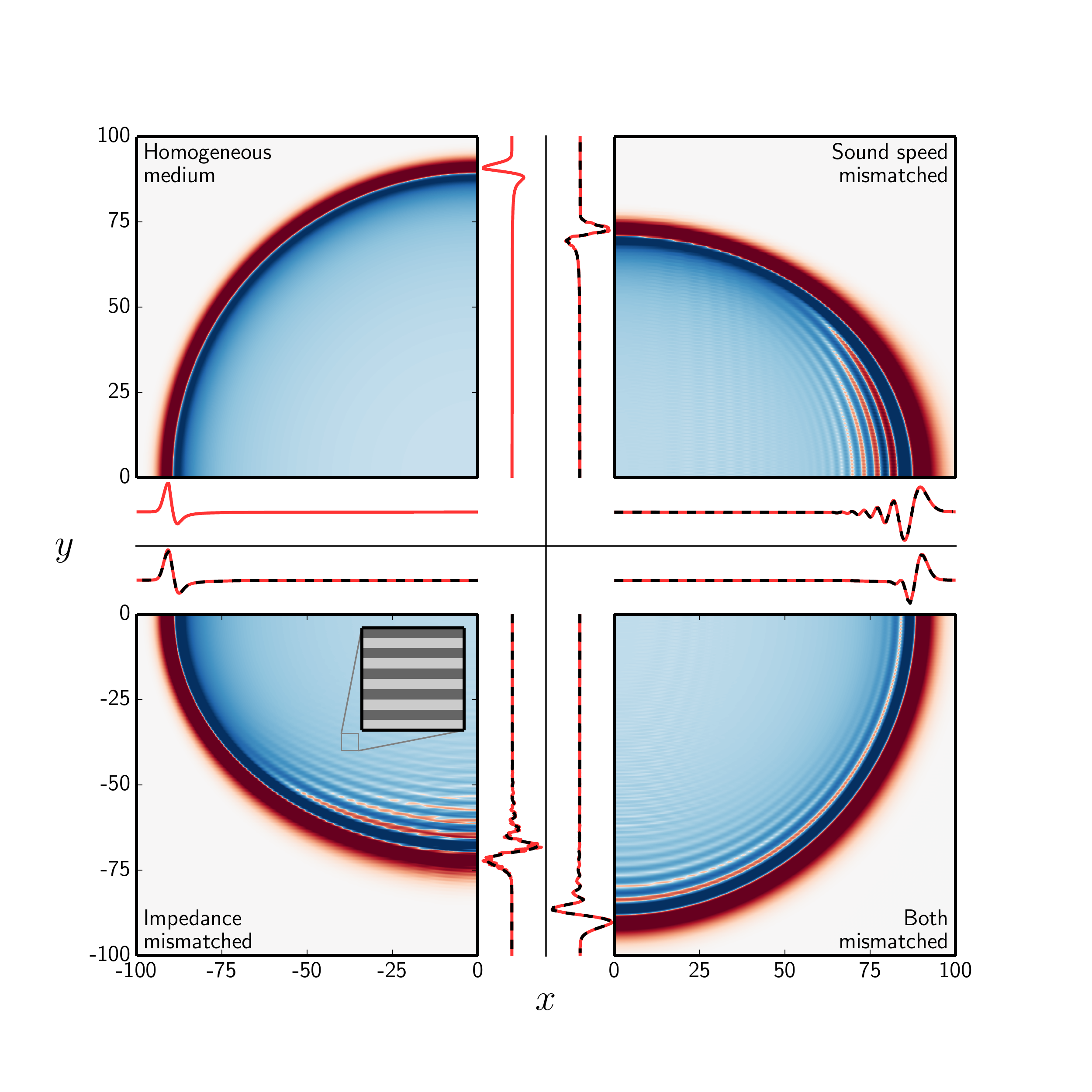}
\caption{Wave propagation in four types of layered periodic media.
The units are scaled to the medium period, and the inset shows the layered
medium structure.  The adjacent line plots show slices of the solution along
the lines $x=0$ (normal) and $y=0$ (transverse).  In the slice plots, the
dotted black lines are approximations based on the effective medium equations
derived in Section \ref{sec:homog}.
Notice how dispersion is evident in the y-direction
when the impedance varies (bottom plots)
and in the x-direction when the sound speed varies (right plots).
\label{fig:quadrants_lin}}
\end{center}
\end{figure}

The qualitative behavior of waves propagating in a periodic medium depends
on whether the sound speed $c(y) = \sqrt{K/\rho}$ and the impedance $Z(y)=\sqrt{K\rho}$
vary in space.  In general both will vary, but special choices can be made
so that either is constant.
Figure \ref{fig:quadrants_lin} shows the typical behavior in each of the
four possible types of media: homogeneous
(top left); constant $Z$ and variable $c$ (top right); constant $c$ and variable $Z$
(bottom left); variable $c$ and $Z$ (bottom right).
Each medium consists of alternating horizontal layers:
\begin{align} \label{discontinuous medium}
 K(y),\rho(y) & = \begin{cases}
    (K_A,\rho_A) \mbox{ if } \left|y-\lfloor y\rfloor-\frac{1}{2}\right|<\frac{1}{4}, \\
    (K_B,\rho_B) \mbox{ if } \left|y-\lfloor y\rfloor-\frac{1}{2}\right|>\frac{1}{4}.
  \end{cases}
\end{align}
In each case, the solution shown corresponds to one quadrant of
the evolution of a an initially Gaussian perturbation:
\begin{align} \label{gaussian}
p_0(x,y) & = e^{-\frac{x^2 + y^2}{2\sigma^2}} & u_0 = v_0 & = 0
\end{align}
with $\sigma=2$.  The line plots show traces of the solution
along the lines $x=0$ (normal) and $y=0$ (transverse).
The solutions are computed using highly-resolved finite volume simulations.

Two important effects are evident.  First, the speed of wave propagation
is {\em anisotropic} in the heterogeneous media; typically, normally-incident waves travel
more slowly.  This is not surprising, given that such waves undergo partial reflection
at each material interface.  This effect is explained well by the lowest
order homogenization theory.
The second effect is that, depending on the nature of the medium,
different components of the wave may develop a dispersive tail.
These {\em dispersive} effects cannot be described by the lowest-order homogenization.

For waves propagating in the normal direction ($\theta=90^\circ$), we observe
that dispersion occurs when the material impedance varies (bottom plots)
and not when the impedance is constant (top plots).  Indeed,
the propagation of a plane wave along the $y$-axis (normal propagation) 
reduces to a one-dimensional problem that is well-studied.
Over long distances, periodic variation in the material impedance 
leads to a dispersive effect -- higher frequencies travel more slowly, due 
to reflection \cite{santosa1991}.  We refer to this as {\em reflective dispersion}.
On the other hand, when the impedance does not
vary, there is no reflection and no effective dispersion
\cite{santosa1991,leveque2003}.
Instead, all wavelengths travel at the harmonic average of the sound speed. 

Next, let us examine the propagation of waves in the transverse direction
($\theta=0^\circ$).  From the left plots, we see that such waves undergo
no dispersion when the material sound speed is constant.  Remarkably, the
right plots show that transverse waves are dispersed when the sound speed 
varies.  In this work we show that diffraction
can play a role similar to that of reflection in periodic media, leading 
again to a dispersive effect in which higher frequencies travel more slowly.
Thus an effective dispersion arises even in materials with constant impedance.
This {\em diffractive dispersion} is an inherently multidimensional effect,
with no one-dimensional analog.
Whereas reflective dispersion depends on variation in the material impedance,
diffractive dispersion depends on variation in the material sound speed.

\begin{figure}
\begin{centering}  
\includegraphics[width=6in]{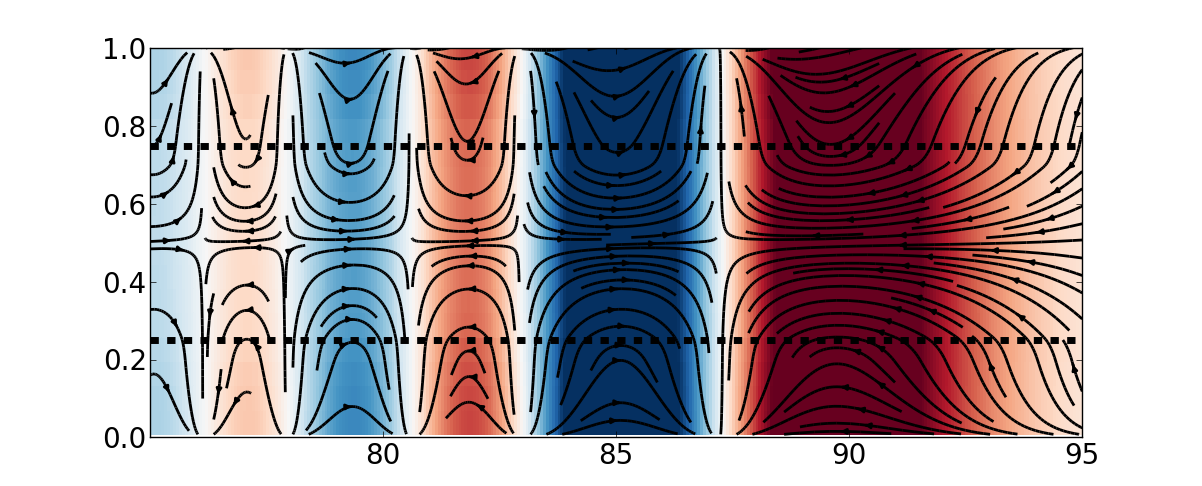}
\caption{Closeup of a transversely-propagating part of the solution for the $c$-dispersive medium
    from the top right quadrant of Figure \ref{fig:quadrants_lin}.
    Streamlines of the velocity field (in black) superimposed on pressure
    (in color).  The wavefront is propagating to the right, but significant diffraction
    is present, as indicated by the vertical velocity components. 
    The dashed lines represent the material interfaces. 
    \label{fig:velocity}}
\end{centering}
\end{figure}

In Figure \ref{fig:velocity}, we plot streamlines of the velocity field superimposed
on a color plot of the pressure, the dashed lines represent the material interfaces.  
The presence of diffraction is evident in the 
velocity streamlines.  It should be noted that the streamlines do not represent 
particle trajectories; they are merely a helpful tool for visualizing the vertical
velocity components created by diffraction.

The dotted black lines in the slice plots are based on pseudospectral solutions
of the effective medium equations \eqref{homog-combined homog system with 4 corrections}, derived in Section \ref{sec:homog} (no corresponding line is shown in the top left quadrant, since the ``effective
medium'' is exact in the homogeneous case).
At the resolution shown in the figure, they are indistinguishable from the
``exact'' finite volume solutions.  The agreement between the effective and variable-coefficient
equations is explored in Section \ref{sec:waves}.

It has been observed that reflective dispersion can, in combination with nonlinear
effects, lead to the formation of solitary waves \cite{leveque2003}.
Diffractive dispersion can also lead to formation of nonlinear solitary waves; 
this is the subject of current research \cite{quezada_diffractons}.

\section{Homogenization\label{sec:homog}}
Homogenization theory can be used to derive an effective PDE for waves in a periodic
medium when the wavelength $\lambda$ is larger than the period of the medium $\Omega$.
The effective PDE is derived through a perturbation expansion, using
$\delta=\Omega/\lambda$ as a small parameter.
We will see that the homogenized PDE depends only on $\Omega$ and not on $\lambda$.
Since the effective PDE has constant coefficients, it can be used to determine
an effective dispersion relation for plane waves in the periodic medium.  

The lowest-order homogenized equation for \eqref{wave-equation} (containing
only terms of $\order(\delta^0)$) is well understood already, but it contains
no dispersive terms and so cannot describe even qualitatively the results shown
in the introduction.
In this section we derive homogenized equations including terms up to $\order(\delta^4)$,
which include dispersive terms.
Additional terms up to $\order(\delta^6)$ are derived for a plane wave
propagating in the $x$-direction. 
Our approach is based on the technique used in \cite{Fish2001} for
one-dimensional wave propagation. 
We will see that some care is required to extend the technique used there to
the case of two-dimensional media of the type shown in Figure \ref{fig:layered}. 
We remark that further difficulties arise when considering periodic
media in which the coefficients depend on both $x$ and $y$; we do not pursue high-order
homogenization for such materials here.

It will be convenient to deal with the wave equation \eqref{wave-equation} in first order form:
\begin{subequations} \label{acoustic equations}
 \begin{align}
   p_t + K(y) (u_x+v_y) & = 0, \\
   \rho(y) u_t + p_x & = 0, \\
   \rho(y) v_t + p_y & = 0.
  \end{align}
\end{subequations}
Here $u,v$ are the velocities in the $x-$ and $y-$coordinate directions respectively.

We start by introducing a fast scale $\yhat=\delta^{-1}y$, and adopting the formalism
that $y$ and $\yhat$ are independent.
The scale defined by $\yhat$ is the scale on which the material properties vary, so we formally replace the
$\Omega$-periodic functions $K(y), \rho(y)$ with $\lambda$-periodic functions
$\Khat(\yhat)$ and $\rhohat(\yhat)$, which are independent of the slow scale $y$.
The dependent variables $p, u, v$ are assumed to vary on both the fast and slow scales,
and are assumed to be periodic in $\hat{y}$ with period $\lambda$.
The idea now is to average over the fast scale $\yhat$ to obtain
constant-coefficient equations involving only $y$.  These equations will not capture
the details of the solution on the fast scale, but will include (to some degree) the influence of
the fast scale on the slow scale.

Using the chain rule we
find that $\partial_y\mapsto\partial_y+\delta^{-1}\partial_{\yhat}$. Therefore, 
system \eqref{acoustic equations} becomes: 
\begin{subequations} \label{homog-scaled system}
 \begin{align}
   \Khat^{-1}(\yhat) p_t + u_x+v_y+\delta^{-1}v_{\yhat} & = 0, \\
   \rhohat (\yhat) u_t + p_x & = 0, \\
   \rhohat (\yhat) v_t + p_y + \delta^{-1} p_{\yhat} & = 0.
  \end{align}
\end{subequations}
For simplicity, from here on, we omit the hats over the coefficient function $K$ and $\rho$, with
the understanding that all material coefficients are $\lambda$-periodic functions
of $\hat{y}$.
Next we assume that $p, u$, and $v$ can be formally written as power series in $\delta$; e.g.
$p(x,y,\yhat,t)=\sum_{i=0}^{\infty}\delta^ip_i(x,y,\yhat,t)$.
Plugging these expansions into \eqref{homog-scaled system} yields
\begin{subequations} \label{homog-expanded system}
 \begin{align}
   K^{-1} \sum_{i=0}^{\infty}\delta^i p_{i,t} + \sum_{i=0}^{\infty}\delta^i u_{i,x}+\sum_{i=0}^{\infty}\delta^i v_{i,y}+\delta^{-1}\sum_{i=0}^{\infty}\delta^i v_{i,{\yhat}} & = 0, \\
   \rho  \sum_{i=0}^{\infty}\delta^i u_{i,t} + \sum_{i=0}^{\infty}\delta^i p_{i,x} & = 0, \\
   \rho  \sum_{i=0}^{\infty}\delta^i v_{i,t} + \sum_{i=0}^{\infty}\delta^i p_{i,y} + \delta^{-1} \sum_{i=0}^{\infty}\delta^i p_{i,{\yhat}} & = 0,
  \end{align}
\end{subequations}
where $(\cdot)_{i,x}$ denotes differentiation of $(\cdot)_i$ with respect to $x$. 
Next we equate terms of the same order in $\delta$;
at each order we apply the averaging operator 
\begin{align} \label{avg}
\avg[\cdot] := \frac{1}{\lambda}\int_0^\lambda(\cdot)d\yhat
\end{align}
to obtain the homogenized leading order system and corrections to it. 
Thus the 
homogenized equations don't depend on the fast scale $\yhat$. 
Note that the averaging operator averages over one period in $y$; i.e., from $0$ to $\Omega$; 
therefore, in $\hat{y}$ it averages from $0$ to $\lambda$; 
i.e., since $y = \delta \yhat =\frac{\Omega}{\lambda}\yhat$, we have
\begin{align} 
\frac{1}{\Omega}\int_0^\Omega(\cdot)dy=\frac{1}{\lambda}\int_0^\lambda(\cdot)d\yhat.
\end{align}
For brevity in longer equations, we will sometimes use a bar instead of
brackets to denote this average (i.e., $\bar{f} = \avg[f]$).

We present the derivation of the first two corrections in detail.
Since the derivation of higher-order terms is similar (but increasingly
tedious), we give the higher-order results without detailed derivations.
Most of the process is mechanical, but for each system we must make an intelligent ansatz 
to obtain an expression for the non-homogenized solution of the corresponding system. 

\subsection{Homogenized $\order(1)$ system}
Taking only $\order(\delta^{-1})$ terms in \eqref{homog-expanded system} gives
\begin{subequations} \label{homog-order d^-1 system}
 \begin{align}
   v_{0,\yhat} & = 0, \\
   p_{0,\yhat} & = 0,
  \end{align}
\end{subequations}
which implies $v_0=:\bar{v}_0(x,y,t)$ and
$p_0=:\bar{p}_0(x,y,t)$ .  Thus the leading order pressure and vertical velocity
are independent of the fast scale $\yhat$. Note that we can't assume that $u_0$
is independent of $\yhat$; we will see in the following sections that it is
not.  This is in contrast to the homogenization of similar systems in 1D where,
to leading order, all dependent variables are independent of the fast scale
\cite{Fish2001}.

Taking only the $\order(1)$ terms in \eqref{homog-expanded system} gives
\begin{subequations} \label{homog-order 1 system}
 \begin{align}
   K^{-1} \bar{p}_{0,t} + u_{0,x}+\bar{v}_{0,y}+v_{1,\yhat} & = 0, \\
   \rho u_{0,t} + \bar{p}_{0,x} & = 0, \label{homog-order 1 system-eqn2} \\
   \rho \bar{v}_{0,t} + \bar{p}_{0,y} + p_{1,\yhat} & = 0.
  \end{align}
\end{subequations}

Next we apply the averaging operator $\avg[\cdot]$ to \eqref{homog-order 1 system}.
This eliminates the terms $v_{1,\hat{y}}$ and $p_{1,\hat{y}}$, which are periodic with
mean zero.
We have no way to determine the average of $\rho u_{0,t}$ because both $\rho$ and $u_0$ depend
on $\yhat$.  We therefore divide
\eqref{homog-order 1 system-eqn2} by $\rho$ and then apply $\avg[\cdot]$, yielding
\begin{subequations} \label{homog-order 1 homog system}
 \begin{align}
   K_h^{-1} \bar{p}_{0,t} + \bar{u}_{0,x}+\bar{v}_{0,y} & = 0, \\
   \rho_h \bar{u}_{0,t} + \bar{p}_{0,x} & = 0, \label{homog-order 1 homog system-eqn2} \\
   \rho_m \bar{v}_{0,t} + \bar{p}_{0,y} & = 0, \label{homog-order 1 homog system-eqn3}
  \end{align}
\end{subequations}
where 
Here and elsewhere the subscripts $m$ and $h$ denote the arithmetic and harmonic average, respectively:
\begin{align*}
\rho_m & := \avg[\rho], & \rho_h & := \avg[\rho^{-1}]^{-1}; \\
K_m & := \avg[K], & K_h & := \avg[K^{-1}]^{-1}.
\end{align*}
We see already that the effective medium is anisotropic, as indicated by the
appearance of these different averages of $\rho$ in 
\eqref{homog-order 1 homog system-eqn2} and \eqref{homog-order 1 homog system-eqn3}.
In particular, we see that plane waves propagating parallel to the $y$-axis 
travel with speed $\sqrt{\frac{K_h}{\rho_m}}$ while
plane waves propagating parallel to the $x$-axis travel with speed
$\sqrt{\frac{K_h}{\rho_h}}$.
We discuss this in more detail in section \ref{sec: effective sound speed}. 

Combining \eqref{homog-order 1 system-eqn2} and \eqref{homog-order 1 homog system-eqn2} yields
\begin{equation} \label{homog-fast changing u0}
 u_0=\frac{\rho_h}{\rho(\yhat)}\bar{u}_0+c,
\end{equation}
where $c$ is a time independent constant. We choose
$\bar{u}_0(x,y,t=0)=\frac{\rho(\hat{y})}{\rho_h}u_0(x,y,\hat{y},t=0)$ so that $c=0$. 
This confirms that $u_0$ varies on the fast scale $\yhat$. More importantly, this
indicates that propagation in $x$ is affected by the heterogeneity in $y$ even
at the macroscopic scale.

Next we obtain expressions for $u_1$, $v_1$ and $p_1$ 
in \eqref{homog-order 1 system}. To do so, we use the following ansatz: 
\begin{subequations}  \label{homog-ansatz for O1}
 \begin{align}
   v_1 & = \bar{v}_1+ A(\yhat) \bar{u}_{0,x}+ B(\yhat)\bar{v}_{0,y}, \\
   p_1 & = \bar{p}_1+ C(\yhat)\bar{p}_{0,y}. \label{homog-ansatz for O1 eqn2}
\end{align}
\end{subequations}
This ansatz is chosen in order to reduce system \eqref{homog-order 1 system} to a system of ODEs.
Substituting the ansatz \eqref{homog-ansatz for O1}, the relation
\eqref{homog-fast changing u0}, and the homogenized leading order system
\eqref{homog-order 1 homog system}
into the the $\order(1)$ system \eqref{homog-order 1 system} we get:
\begin{align}
	(A_{\hat{y}}-K_hK^{-1}+\rho_h\rho^{-1})\bar{u}_{0,x}+(B_{\hat{y}}-K_hK^{-1}+1)\bar{v}_{0,x} & = 0, \\
	(C_{\yhat} - \rho\rho_m^{-1}+1)\bar{p}_{0,y} & = 0.
\end{align}
Based on this, it is convenient to choose $A(\hat{y})$, $B(\hat{y})$ and $C(\hat{y})$ to satisfy
\begin{subequations} \label{homog-ODEs for A, B and C}
\begin{align}
 A_{\yhat} - K^{-1}K_h + \rho^{-1}\rho_h & = 0, \\
 B_{\yhat} - K^{-1}K_h-1 & = 0, \\
 C_{\yhat}  - \rho\rho_m^{-1}+1& =0.
\end{align}
\end{subequations}
Equations \eqref{homog-ODEs for A, B and C} represent boundary value ODEs with 
the normalization conditions that $\avg[ A] =\avg[ B] =\avg[ C] =0$.
Note that $\avg[ A_{\yhat}]=\avg[ B_{\yhat}]=\avg[ C_{\yhat}]=0$,
which implies that $A$, $B$ and $C$ are $\lambda$-periodic. 
To solve these boundary value problems we must specify the material functions $\rho,K$. 
In Appendix \ref{sec: fast variable functions} we show the fast-variable functions 
$A$, $B$ and $C$ for a layered medium.
It is convenient to introduce the following linear operators (see \cite[p. 1554]{leveque2003}):
\begin{subequations}
\begin{align}
    \{a\}(\yhat) & = a(\yhat) - \avg[a(\yhat)] \\
    \intavg[a] (\yhat) & = \int_s^{\yhat} \{a\}(\xi)d\xi, & \text{ where $s$ is chosen such that } \avg[\intavg[a] (\yhat)] &= 0.
\end{align}
\end{subequations}
Then we have
\begin{subequations}
\begin{align} \label{A, B and C}
    A & = \intavg[K^{-1}K_h-\rho^{-1}\rho_h] \\
    B & = \intavg[K^{-1}K_h],\\
    C & = \intavg[\rho\rho_m^{-1}].
\end{align}
\end{subequations}

\subsection{Derivation of $\order(\delta)$ system}
Taking only the $\order(\delta)$ terms in \eqref{homog-expanded system} gives
\begin{subequations} \label{homog-order delta system}
 \begin{align}
   K^{-1} p_{1,t} + u_{1,x}+v_{1,y}+v_{2,\yhat} & = 0, \\
   \rho u_{1,t} + p_{1,x} & = 0, \label{homog-order delta system-eqn2} \\
   \rho v_{1,t} + p_{1,y} + p_{2,\yhat} & = 0.
  \end{align}
\end{subequations}
Substituting the ansatz for $v_1$ and $p_1$ from \eqref{homog-ansatz for O1} into 
\eqref{homog-order delta system} and averaging gives
\begin{subequations}
\begin{align}
 K_h^{-1}\bar{p}_{1,t}+\bar{u}_{1,x}+\bar{v}_{1,y} & = -\avg[ K^{-1}C] \bar{p}_{0,yt}, \\
 \rho_h  \bar{u}_{1,t}+\bar{p}_{1,x} & = -\rho_h\avg[ \rho^{-1}C] \bar{p}_{0,xy},\\
 \rho_m \bar{v}_{1,t}+\bar{p}_{1,y} & = -\avg[ \rho A] \bar{u}_{0,xt}-\avg[ \rho B] \bar{v}_{0,yt},
 \end{align}
\end{subequations}
where $\bar{u}_1:=\avg[u_1]$ and similarly for $\bar{v}_1$ and $\bar{p}_1$.
For any piecewise-constant functions $f,g$ it can be shown that $\avg[ f\llbracket g \rrbracket ]=0$
\cite{leveque2003}.
Thus, for the piecewise-constant materials that we consider in Section \ref{sec:waves}, we have
$\avg[ K^{-1}C] =\avg[ \rho^{-1}C] =\avg[ \rho A] =\avg[ \rho B] =0$. 
These averages also vanish for the sinusoidal materials we will consider.
For more general materials, these terms may be non-zero, but in the following
we assume they vanish.  Then we obtain:
\begin{subequations} \label{homog-order delta homog system}
\begin{align}
 K_h^{-1}\bar{p}_{1,t}+\bar{u}_{1,x}+\bar{v}_{1,y} & = 0, \\
 \rho_h  \bar{u}_{1,t}+\bar{p}_{1,x} & = 0,\\
 \rho_m \bar{v}_{1,t}+\bar{p}_{1,y} & = 0.
 \end{align}
\end{subequations}
Since the boundary conditions are imposed in the leading order homogenized 
system \eqref{homog-order 1 homog system}, system \eqref{homog-order delta homog system} 
should be solved with homogeneous Dirichlet boundary conditions; therefore, its
solution vanishes:
\begin{equation}
 \bar{u}_1= \bar{v}_1 = \bar{p}_1 = 0.
\end{equation}

Taking $\bar{u}_1=\bar{v}_1=\bar{p}_1=0$, we make the following ansatz 
for the $\order(\delta)$ solutions $v_2$ and $p_2$:
\begin{align}\label{homog-ansatz for O delta}
 v_2 & = \bar{v}_2+D(\yhat)\bar{u}_{0,xy}+E(\yhat)\bar{v}_{0,yy},\nonumber \\
 p_2 & = \bar{p}_2+F(\yhat)\bar{p}_{0,yy}+H(\yhat)\bar{p}_{0,xx},
\end{align}
which is chosen in order to reduce system \eqref{homog-order delta system} 
to a system of ODEs. 
From \eqref{homog-order delta system-eqn2} we get $u_{1,t}=-\rho^{-1}p_{1,x}$. 
The ansatz for $p_1$ from \eqref{homog-ansatz for O1 eqn2} gives
$u_{1,t}=-\rho^{-1}C(\bar{p}_{0,x})_y$ 
and using the homogenized leading order equation \eqref{homog-order 1 homog system-eqn2} we
get $u_{1,t}=\rho^{-1}\rho_hC(\bar{u}_{0,y})_t$. Finally we get an expression for the non-homogenized solution $u_1$:
\begin{equation} \label{homog-fast changing}
 u_1=\rho^{-1}\rho_hC\bar{u}_{0,y}.
\end{equation}

Next we substitute the ansatz for $v_1$ and $p_1$ from \eqref{homog-ansatz for O1},
the ansatz for $v_2$ and $p_2$ from \eqref{homog-ansatz for O delta} and
the non-homogenized solution $u_1$ from \eqref{homog-fast changing} into system
\eqref{homog-order delta system}. 
Then we substitute the leading order homogenized system \eqref{homog-order 1 homog system}
and the coefficients \eqref{A, B and C}. Finally, in the resulting expression, 
we make the fast variable coefficients to vanish to obtain
\begin{subequations}\label{D, E, F and H}
\begin{align}
    D & = \intavg[ K^{-1}K_hC-\rho^{-1}\rho_hC-A],\\
    E & = \intavg[ K^{-1}K_hC-B],\\
    F & = \intavg[ \rho\rho_m^{-1}B-C],\\
    H & = \intavg[ \rho\rho_h^{-1} A].
\end{align}
\end{subequations}

\subsection{Derivation of $\order(\delta^2)$ system}
From \eqref{homog-expanded system} we take $\order(\delta^2)$ terms to get:
\begin{subequations} \label{homog-order delta2 system}
 \begin{align}
   K^{-1} p_{2,t} + u_{2,x}+v_{2,y}+v_{3,\yhat} & = 0, \\
   \rho u_{2,t} + p_{2,x} & = 0, \\
   \rho v_{2,t} + p_{2,y} + p_{3,\yhat} & = 0.
  \end{align}
\end{subequations}
Substituting the ansatz for $v_2$ and $p_2$ from \eqref{homog-ansatz for O delta}
into \eqref{homog-order delta2 system} and averaging yields
\begin{subequations} \label{homog-order delta2 homog system}
\begin{align} 
 K_h^{-1}\bar{p}_{2,t}+\bar{u}_{2,x}+\bar{v}_{2,y} & = \alpha_1 (\bar{u}_{0,xyy}+\bar{v}_{0,yyy})+\alpha_2 (\bar{u}_{0,xxx}+\bar{v}_{0,xxy}), \\ 
 \rho_h\bar{u}_{2,t}+\bar{p}_{2,x} & = \beta_1 \bar{p}_{0,xyy}+\beta_2 \bar{p}_{0,xxx}, \\
 \rho_m\bar{v}_{2,t}+\bar{p}_{2,y} & = \gamma_1 \bar{p}_{0,yyy}+\gamma_2 \bar{p}_{0,xxy},
\end{align}
\end{subequations}
where $\bar{u}_2:=\avg[u_2]$ and similarly for $\bar{v}_2$ and $\bar{p}_2$. Expressions for the coefficients $\alpha, \beta, \gamma$ appear in
appendix \ref{sec: coefficients of homogenized equations}.
 
\subsection{Higher order corrections}
Following similar, but more involved steps, we find the $\order(\delta^3)$ and $\order(\delta^4)$ corrections.

\subsubsection{$\order(\delta^3)$ homogenized correction}

The third homogenized correction is:
\begin{subequations} \label{homog-order delta3 homog system}
\begin{align} 
 K_h^{-1}\bar{p}_{3,t}+\bar{u}_{3,x}+\bar{v}_{3,y} & = -\avg[ K^{-1}C ] \bar{p}_{2,yt} 
   									    -\avg[ K^{-1}N] \bar{p}_{0,yyyt}
									    -\avg[ K^{-1}P] \bar{p}_{0,xxyt}, \\
 \rho_h\bar{u}_{3,t}+\bar{p}_{3,x} & = -\rho_h\avg[ \rho^{-1}C] \bar{p}_{2,xy}  
				        -\rho_h \avg[ \rho^{-1}N] \bar{p}_{0,xyyy} 
					-\rho_h \avg[ \rho^{-1}P] \bar{p}_{0,xxxy}, \\
 \rho_m\bar{v}_{3,t}+\bar{p}_{3,y} & = -\avg[ \rho A] \bar{u}_{2,xt} 
					- \avg[ \rho B] \bar{v}_{2,yt} \\
			 	     & - \avg[ \rho I] \bar{u}_{0,xyyt} 
					- \avg[ \rho J] \bar{u}_{0,xxxt}
					- \avg[ \rho L] \bar{v}_{0,xxyt}
					- \avg[ \rho M] \bar{v}_{0,yyyt}, \nonumber
\end{align}
\end{subequations}
where $\bar{u}_3:=\avg[u_3]$ and similarly for $\bar{v}_3$ and $\bar{p}_3$. The fast-variable functions $I(\yhat)$, $J(\yhat)$, $L(\yhat)$,
$M(\yhat)$, $N(\yhat)$, $P(\yhat)$ are solutions of the BVPs
\begin{subequations}\label{I, J, L, M, P and N}
\begin{align}
I & = \llbracket K^{-1}K_h\left(F-K_h\avg[ K^{-1}F ] \right) -\rho^{-1}\rho_h\left(F-\rho_h\avg[ \rho^{-1}F ] \right) -D\rrbracket,\\
J& = \llbracket K^{-1}K_h\left(H-K_h\avg[ K^{-1}H ] \right) 
        -\rho^{-1}\rho_h\left(H-\rho_h\avg[ \rho^{-1}H ] \right)\rrbracket,\\
L& = \llbracket K^{-1}K_h\left(H-K_h\avg[ K^{-1}H ] \right)\rrbracket, \\ 
M& = \llbracket K^{-1}K_h\left(F-K_h\avg[ K^{-1}F ] \right)-E\rrbracket, \\ 
N& = \llbracket \rho \rho_m^{-1}\left(E-\rho_m^{-1}\avg[ \rho E ] \right)-F\rrbracket, \\ 
P& = \llbracket \rho \rho_h^{-1}\left(D-\rho_m^{-1}\avg[ \rho D ] \right)-H\rrbracket.
\end{align}
\end{subequations}
For the two types of media considered in this work all coefficients on the right hand side 
of \eqref{homog-order delta3 homog system} vanish. Since the boundary conditions
are fulfilled by the leading order homogenized system \eqref{homog-order 1 homog system}, 
the third homogenized correction vanishes; i.e., 
\begin{equation}
 \bar{u}_3= \bar{v}_3 = \bar{p}_3 = 0.
\end{equation}

\subsubsection{$\order(\delta^4)$ homogenized correction}

The fourth correction is given by: 
\begin{subequations} \label{homog-order delta4 homog system}
\begin{align} 
 K_h^{-1}\bar{p}_{4,t}+\bar{u}_{4,x}+\bar{v}_{4,y} & = \alpha_1(\bar{u}_{2,xyy}+\bar{v}_{2,yyy})+\alpha_2(\bar{u}_{2,xxx}+\bar{v}_{2,xxy}) \nonumber \\ 
 & + \alpha_3 (\bar{u}_{0,xyyyy}+\bar{v}_{0,yyyyy} ) + \alpha_4 (\bar{u}_{0,xxxxx}+\bar{v}_{0,xxxxy} ) \nonumber \\
 & + \alpha_5 (\bar{u}_{0,xxxyy}+\bar{v}_{0,xxyyy} ), \\
 \rho_h\bar{u}_{4,t}+\bar{p}_{4,x} & = \beta_1 \bar{p}_{2,xyy} + \beta_2 \bar{p}_{2,xxx} \nonumber \\
 & + \beta_3 \bar{p}_{0,xyyyy} + \beta_4 \bar{p}_{0,xxxxx} + \beta_5 \bar{p}_{0,xxxyy}, \\
 \rho_m\bar{v}_{4,t}+\bar{p}_{4,y} & = \gamma_1\bar{p}_{2,yyy} + \gamma_2 \bar{p}_{2,xxy} \nonumber \\
 & + \gamma_3 \bar{p}_{0,yyyyy} + \gamma_4 \bar{p}_{0,xxyyy} + \gamma_5 \bar{p}_{0,xxxxy},
\end{align}
\end{subequations}
where $\bar{u}_4:=\avg[u_4]$ and similarly for $\bar{v}_4$ and $\bar{p}_4$.
Expressions for the coefficients $\alpha, \beta, \gamma$ are given in
appendix \ref{sec: coefficients of homogenized equations} 

\subsection{Combined homogenized equations}
Once we have the homogenized leading order system and the homogenized
corrections we can combine them into a single system. This is done by taking 
$\bar{p}:=\avg[ p_0+\delta p_1 + \dots + \delta^4 p_4 ]=\bar{p}_0+\delta \bar{p}_1 + \dots + \delta^4 \bar{p}_4$
and similarly for $\bar{u}$ and $\bar{v}$. 
Combining homogenized systems \eqref{homog-order 1 homog system}, 
\eqref{homog-order delta homog system}, \eqref{homog-order delta2 homog system}, \eqref{homog-order delta3 homog system} 
and \eqref{homog-order delta4 homog system} we obtain: 
\begin{subequations} \label{homog-combined homog system with 4 corrections}
\begin{align}\label{24a}  
 K_h^{-1}\bar{p}_t+\bar{u}_x+\bar{v}_y & = \delta^2\left[ \alpha_1(\bar{u}_{xyy}+\bar{v}_{yyy})+\alpha_2(\bar{u}_{xxx}+\bar{v}_{xxy})\right] \nonumber \\ 
 & \quad + \delta^4 \left[\alpha_3 (\bar{u}_{xyyyy}+\bar{v}_{yyyyy} ) + \alpha_4 (\bar{u}_{xxxxx}+\bar{v}_{xxxxy} )\right] \\
 & \quad + \delta^4 \left[\alpha_5 (\bar{u}_{xxxyy}+\bar{v}_{xxyyy} )\right], \nonumber \\
 \rho_h \bar{u}_{t}+\bar{p}_x & = \delta^2\left[\beta_1 \bar{p}_{xyy} + \beta_2 \bar{p}_{xxx}\right] + \delta^4\left[\beta_3 \bar{p}_{xyyyy} + \beta_4 \bar{p}_{xxxxx} + \beta_5 \bar{p}_{xxxyy}\right], \label{24b} \\
 \rho_m \bar{v}_t+\bar{p}_y & = \delta^2\left[\gamma_1\bar{p}_{yyy} + \gamma_2 \bar{p}_{xxy}\right] + \delta^4\left[\gamma_3 \bar{p}_{yyyyy} + \gamma_4 \bar{p}_{xxyyy} + \gamma_5 \bar{p}_{xxxxy}\right], \label{24c}
\end{align}
\end{subequations}
where expressions for the coefficients are given in appendix 
\ref{sec: coefficients of homogenized equations}. 
Unlike the lowest-order homogenized equation
\eqref{homog-order 1 system}, this system is dispersive and 
the dispersion depends on the direction of propagation, 
see section \ref{sec:dispersion}.

In general, each coefficient of a 
$\delta^n=\frac{\Omega^n}{\lambda^n}$ term in
\eqref{homog-combined homog system with 4 corrections}
contains a matching factor $\lambda^n$ (contained within the 
$\alpha^\prime s$ and $\beta^\prime s$; see appendix
\ref{sec: coefficients for piecewise constant medium}
for an example of the coefficients of the first order solution 
for a layered medium).
As a result, each 
term on the right hand side of 
\eqref{homog-combined homog system with 4 corrections}
is proportional to $\Omega^n$ 
(and independent of $\lambda$).
This explains the observation in \cite{leveque2003} that the 
homogenized equations are valid for any choice of the material
period $\Omega$.  
Nevertheless, \eqref{homog-combined homog system with 4 corrections} is only valid for small 
$\delta$; i.e., for relatively long wavelengths $\lambda>\Omega$. 
In all numerical simulations in this work, we take $\Omega=1$.


\section{Effective dispersion relations\label{sec:dispersion}}
Equation \eqref{homog-combined homog system with 4 corrections} is a 
linear system of PDEs with constant coefficients.  Hence its solutions
can be completely described by the dispersion relation,
which relates frequency and wavenumber for a plane wave:
\begin{equation}\label{plane wave ansatz}
 \bar{p}(x,y,t)=\bar{p}_0e^{i\left({\bf k}\cdot{\bf x} - \omega t\right)}.
\end{equation}
Here $\bar{p}_0$ is the amplitude, $\omega$ is the angular frequency and 
${\bf k}$ is the wave vector. Let ${\bf k}=k(k_x, k_y)$ with 
$k_x = \cos\theta, k_y=\sin\theta$, where $\theta$ is the direction of propagation
(see Figure \ref{fig:layered}).
Because \eqref{homog-combined homog system with 4 corrections} describes
a medium that is anisotropic and dispersive, the speed of propagation
of a plane wave depends on both the angle $\theta$ and the wavenumber magnitude $k$.

We can combine \eqref{homog-combined homog system with 4 corrections} into
a single second-order equation by differentiating \eqref{24a} with respect to
$t$, differentiating \eqref{24b} and \eqref{24c} with respect to $x$ and $y$ 
respectively, and equating mixed partial derivatives.  By
substituting \eqref{plane wave ansatz} into the result, we obtain the
effective dispersion relation, up to $\order(\delta^4)$:
\begin{equation} \label{effective dispersion}
\begin{split}
 \omega^2  &= \frac{K_h}{\rho_h\rho_m}k^2 \left(k_x^2\rho_m+k_y^2\rho_h\right)
                 + \delta^2\frac{K_h}{\rho_h\rho_m}k^4 \left[k_x^2\rho_m\left(\left(\alpha_2+\beta_2\right)k_x^2+\left(\alpha_1+\beta_1\right)k_y^2\right) \right. \\
                & \left. \vphantom{} +k_y^2\rho_h\left(\left(\alpha_2+\gamma_2\right)k_x^2+\left(\alpha_1+\gamma_1\right)k_y^2\right)\right] 
                  + \delta^4\frac{K_h}{\rho_h\rho_m}k^6 [ k_x^2\rho_m\left(\left(\alpha_4+\alpha_2\beta_2-\beta_4\right)k_x^4 \right. \\
                & \left. \vphantom{} -\left(\alpha_5-\alpha_2\beta_1-\alpha_1\beta_2+\beta_5\right)k_x^2k_y^2+\left(\alpha_3+\alpha_1\beta_1-\beta_3\right)k_y^4\right) \\
                & + k_y^2\rho_h\left(\left(\alpha_4+\alpha_2\gamma_2-\gamma_5\right)k_x^4-\left(\alpha_5-\alpha_2\gamma_1-\alpha_1\gamma_2+\gamma_4\right)k_x^2k_y^2
                +\left(\alpha_3+\alpha_1\gamma_1-\gamma_3\right)k_y^4\right) ].
\end{split} 
\end{equation}
It is important to keep in mind that although \eqref{effective dispersion}
is accurate to $\order(\delta^4)$, it is obtained by homogenization and
is not expected to be valid for wavelengths shorter than $\Omega$, the medium period.
This is because $\delta=\frac{\Omega}{\lambda}$ is assumed to be small; in particular, $\delta<1$.

\subsection{Effective sound speed} \label{sec: effective sound speed}
Taking only $\order(1)$ terms in \eqref{effective dispersion} we obtain the
effective sound speed: 
\begin{equation} \label{sound speed}
 \ceff = \omega/k = \sqrt{\frac{K_h}{\rho_h} k_x^2 + \frac{K_h}{\rho_m} k_y^2}.
\end{equation}
The effective sound speed, which indicates the speed of very long wavelength 
perturbations, depends on the direction of propagation. 
For normally incident waves ($\theta = 90^\circ$), we have $\ceff = \sqrt{K_h/\rho_m}$,
which is the effective sound speed in a 1D layered medium \cite{santosa1991}.
For transverse waves ($\theta=0^\circ$), we have $\ceff = \sqrt{K_h/\rho_h}$. 
Since the harmonic average is less than or equal to the arithmetic average,
long-wavelength normal waves never travel faster than their transverse wave
counterparts.  This is intuitively reasonable since transverse propagating
waves undergo no reflection.

\begin{figure}
\begin{centering}
  \includegraphics[scale=0.5]{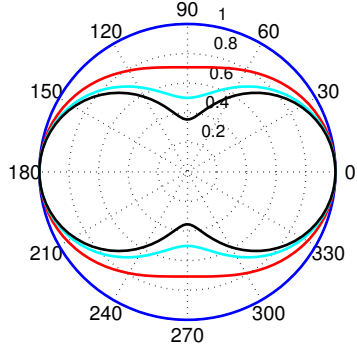}
\par\end{centering}
\caption{Polar plot of the effective sound speed \eqref{sound speed} with $K_h=1$, 
$\rho_h=1$ and $\rho_m=1$ (blue), $\rho_m=2$ (red), $\rho_m=4$ (cyan), 
$\rho_m=8$ (black).}
\label{fig: sound speed}
\end{figure}

\begin{figure}
\begin{centering}
  \includegraphics[scale=0.425]{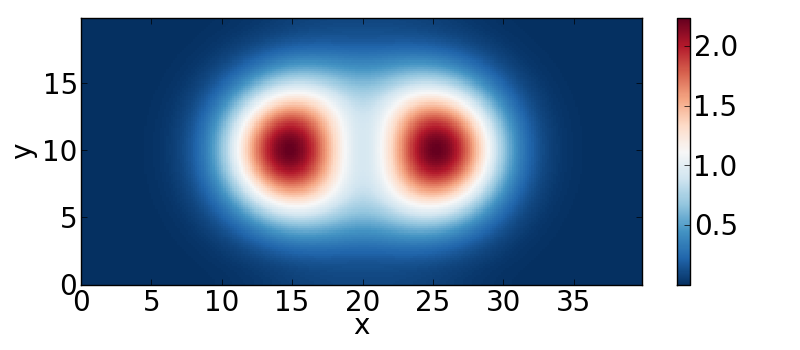}
  \includegraphics[scale=0.425]{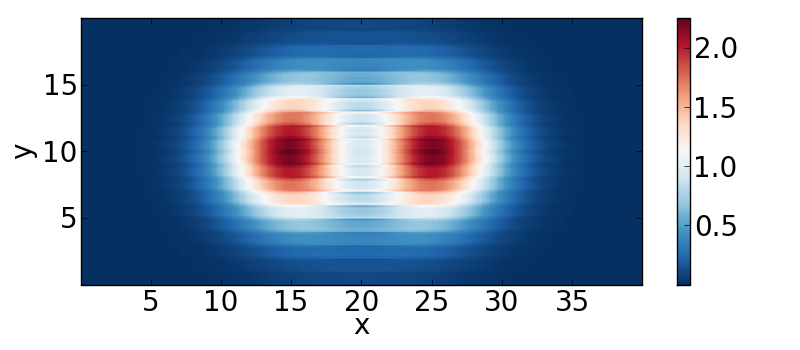}
\par\end{centering}
\caption{Homogenized leading order solution $\bar{p}$ (left) and finite volume solution $p$ (right)
at $t=5$ in a layered medium with $\rho_A=8+\sqrt{56}$, $\rho_B=8/ \rho_A$ and $K_A=K_B=1$ 
(which gives $K_h=\rho_h=1$ and $\rho_m=8$).
This corresponds to the effective sound speed distribution shown in Figure \ref{fig: sound speed} (black line).  Note that the horizontal lines in the right figure are the effect of the
medium structure, and not of numerical resolution (which is much finer).}
\label{fig: sound speed HOM-FV solution}
\end{figure}

In Figure \ref{fig: sound speed} we plot $\ceff$ as a function of $\theta$
for material parameters $K_h=1$, $\rho_h=1$ and 
different values of $\rho_m$.  When $\rho_h=\rho_m$, the sound speed is the same in all 
directions so we obtain the blue line in Figure \ref{fig: sound speed}. 
As $\rho_m$ increases (corresponding to strong impedance variation and thus more reflection),
the effective speed in $y$ decreases and we obtain the red, cyan and black
lines in Figure \ref{fig: sound speed}.  In Figure \ref{fig: sound speed HOM-FV solution} 
we take the initial condition
\begin{align}
\bar{p}_0(x,y) & = 5 e^{-\frac{(x-20)^2 + (y-10)^2}{2\sigma^2}} & \bar{u}_0 = \bar{v}_0 & = 0
\end{align}
with $\sigma^2 = 10$
and a piecewise-constant medium \eqref{discontinuous medium}
with $\rho_A=8+\sqrt{56}$, 
$\rho_B=8/ \rho_A$ and $K_A=K_B=1$ which gives $K_h=\rho_h=1$ and $\rho_m=8$, 
corresponding to the black line in Figure \ref{fig: sound speed}. 
We show the homogenized leading order pressure (left) and the 
finite volume pressure (right) at $t=5$, which correspond to the solution of 
\eqref{homog-order 1 homog system} and \eqref{wave-equation} respectively. 
The predicted anisotropic behavior is observed in both solutions.

\subsection{Normally and transversely incident waves\label{sec:directions}}
The dispersion relation \eqref{effective dispersion} has some special properties for waves
that are aligned with the coordinate axes.  Recall the definitions
of normal and transverse wave propagation illustrated in Figure \ref{fig:layered}.

A normally-incident plane wave corresponds to initial data that is constant in $x$.
For such waves, \eqref{homog-combined homog system with 4 corrections} reduces to
a one-dimensional equation:
\begin{subequations} \label{homog-z-dispersion homog system with 4 corrections}
\begin{align} 
 K_h^{-1}\bar{p}_t+\bar{v}_y & = \delta^2 \alpha_1 \bar{v}_{yyy} + \delta^4 \alpha_3 \bar{v}_{yyyyy}, \\
 \rho_m \bar{v}_{t}+\bar{p}_y & = \delta^2 \gamma_1 \bar{p}_{yyy} + \delta^4 \gamma_3 \bar{p}_{yyyyy}.
\end{align}
\end{subequations}
This system was obtained previously in \cite{santosa1991,Fish2001,leveque2003}.
For a piecewise-constant medium, the coefficients on the right hand side are
all proportional to the difference of squared impedance; for instance
\begin{align}
 \alpha_1 & = - \left(Z_A^2-Z_B^2\right) \frac{\left(K_A-K_B\right)}{192K_m^2 \rho_m}\lambda^2,
\end{align}
where $Z_A=\sqrt{K_A\rho_A}$ and $Z_A=\sqrt{K_B\rho_B}$.

In general all of the dispersive terms in 
\eqref{homog-z-dispersion homog system with 4 corrections}
vanish when the impedance is constant (see Appendix 
\ref{sec: coefficients of homogenized equations}).
This is because no reflection occurs in such media.
This explains why the normally-propagating waves in Figure \ref{fig:quadrants_lin}
are dispersed in the lower two plots (with variable $Z$) but not in the
upper two plots (with constant $Z$).

Transversely-incident plane waves correspond to initial data that is constant in $y$.
For such waves, \eqref{homog-combined homog system with 4 corrections} simplifies to
\begin{subequations} \label{homog-c-dispersion homog system with 6 corrections}
\begin{align} 
 K_h^{-1}\bar{p}_t+\bar{u}_x & = \delta^2 \alpha_2 \bar{u}_{xxx} + \delta^4 \alpha_4 \bar{u}_{xxxxx} + \delta^6\alpha_6 \bar{u}_{xxxxxxx}, \\
 \rho_h \bar{u}_{t}+\bar{p}_x & = \delta^2 \beta_2 \bar{p}_{xxx} + \delta^4 \beta_4 \bar{p}_{xxxxx} + \delta^6 \beta_6 \bar{p}_{xxxxxxx}.
\end{align}
\end{subequations}
Here we have included additional 6th-order corrections, because this
case will be of particular interest in what follows.
For a piecewise-constant medium, all coefficients on the right hand side are proportional to the
difference of squared sound speeds; for instance
\begin{align}
 \alpha_2 & = -\left(c_A^2-c_B^2\right)\frac{\left(K_A-K_B\right)}{192K_m^2\rho_m^{-1}}\lambda^2,
\end{align}
where $c_A=\sqrt{\frac{K_A}{\rho_A}}$ and $c_B=\sqrt{\frac{K_B}{\rho_B}}$.

In general, all of the dispersive terms in 
\eqref{homog-c-dispersion homog system with 6 corrections} vanish when the
sound speed is constant
(see Appendix \ref{sec: coefficients of homogenized equations}).
This is because no diffraction occurs for transversely-incident plane waves in such media.
This explains why the transversely-propagating waves in Figure \ref{fig:quadrants_lin}
are dispersed in the right two plots (with variable $c$) but not in the
left two plots (with constant $c$).

Thus we see that, with respect to dispersion, the role played by impedance for
normal wave propagation corresponds to the role played by sound speed for
transverse wave propagation, and vice-versa.

\section{Comparison of effective medium and variable-coefficient medium solutions\label{sec:waves}}
In this section, we compare in greater detail the solutions obtained from the homogenized effective
medium equations \eqref{homog-combined homog system with 4 corrections} and from the variable-coefficient equation \eqref{wave-equation}.
Solutions of the homogenized equations are obtained
using a Fourier spectral collocation method in space and 4th order Runge-Kutta
integration in time \cite{trefethen2000spectral}.  Solutions of the
variable-coefficient problem are obtained using PyClaw \cite{pyclaw-sisc}, including
numerical methods described in \cite{leveque1997,Ketcheson2011,quezada2013}.
The code used to produce these results is available from
\url{http://github.com/ketch/effective_dispersion_RR}.

\subsection{Propagation of a plane-wave perturbation}
In this section we show the effect of $c$-dispersion on an initial
perturbation that varies only in $x$:
\begin{align*} 
\bar{p}_0(x,y) & = 10 e^{-\frac{x^2}{10}} & \bar{u}_0 = \bar{v}_0 & = 0.
\end{align*}

The homogenized equations used are given by \eqref{homog-c-dispersion homog system with 6 corrections}. 
To demonstrate that the homogenized equations are valid for more
general media, we consider both the piecewise constant medium \eqref{discontinuous medium} 
and a smoothly varying medium
\begin{subequations} \label{sinusoidal medium}
\begin{align}
  K(y) & = \frac{K_A+K_B}{2} +
                \frac{K_A-K_B}{2} \sin\left(2\pi y\right), \\
  \rho(y) & = \frac{1}{K(y)},
\end{align}
\end{subequations}
with the material parameters
\begin{align}
    K_A & = \frac{1}{\rho_A} = 5/8, & K_B & = \frac{1}{\rho_B} = 5/2.
\end{align}
For the sinusoidally-varying medium we numerically solve 
the BVPs defining the homogenization coefficients.
Mathematica and Matlab files to solve for these expressions can be found at 
\url{http://github.com/ketch/effective_dispersion_RR}.

\begin{figure}
\begin{centering}
 \subfloat[Layered medium\label{fig: c-disp disc FV vs Hom}]{
  \includegraphics[scale=0.3]{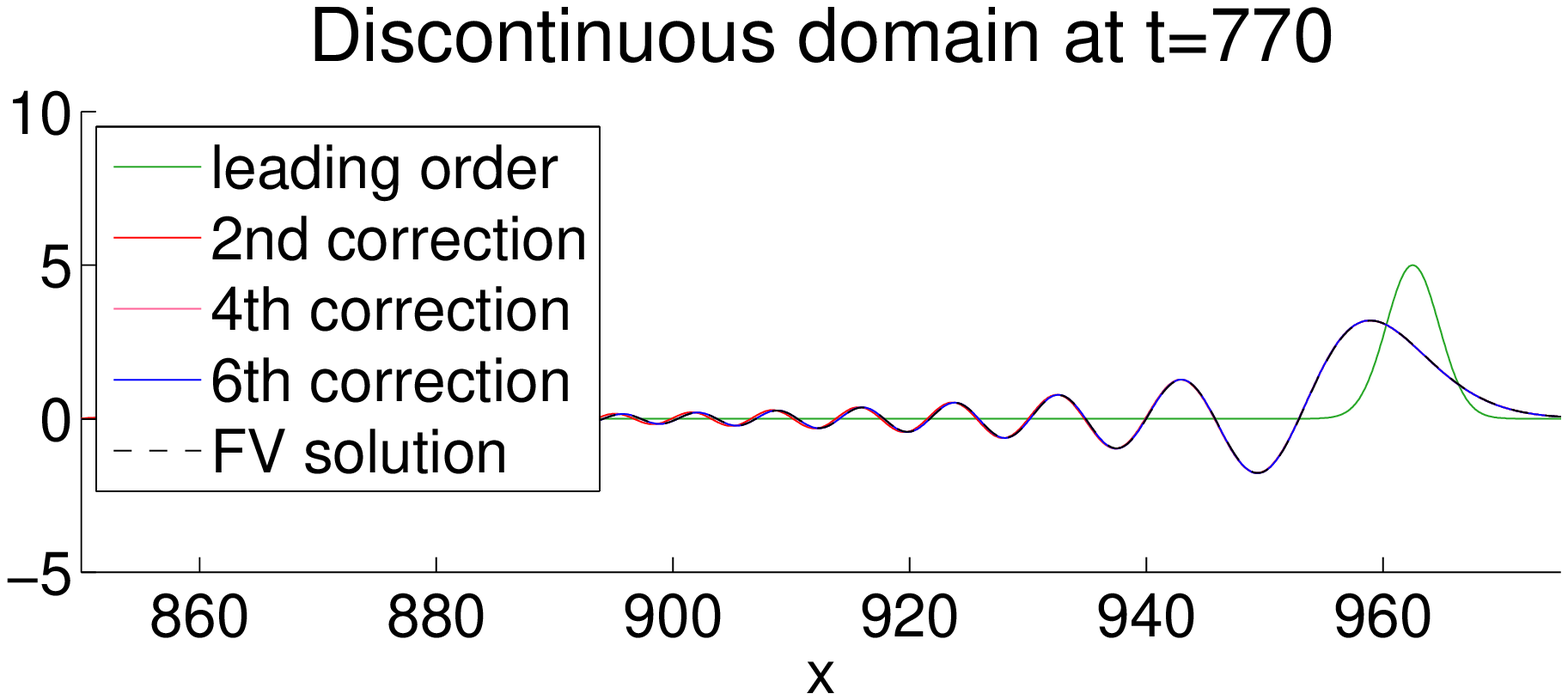}
  \includegraphics[scale=0.3]{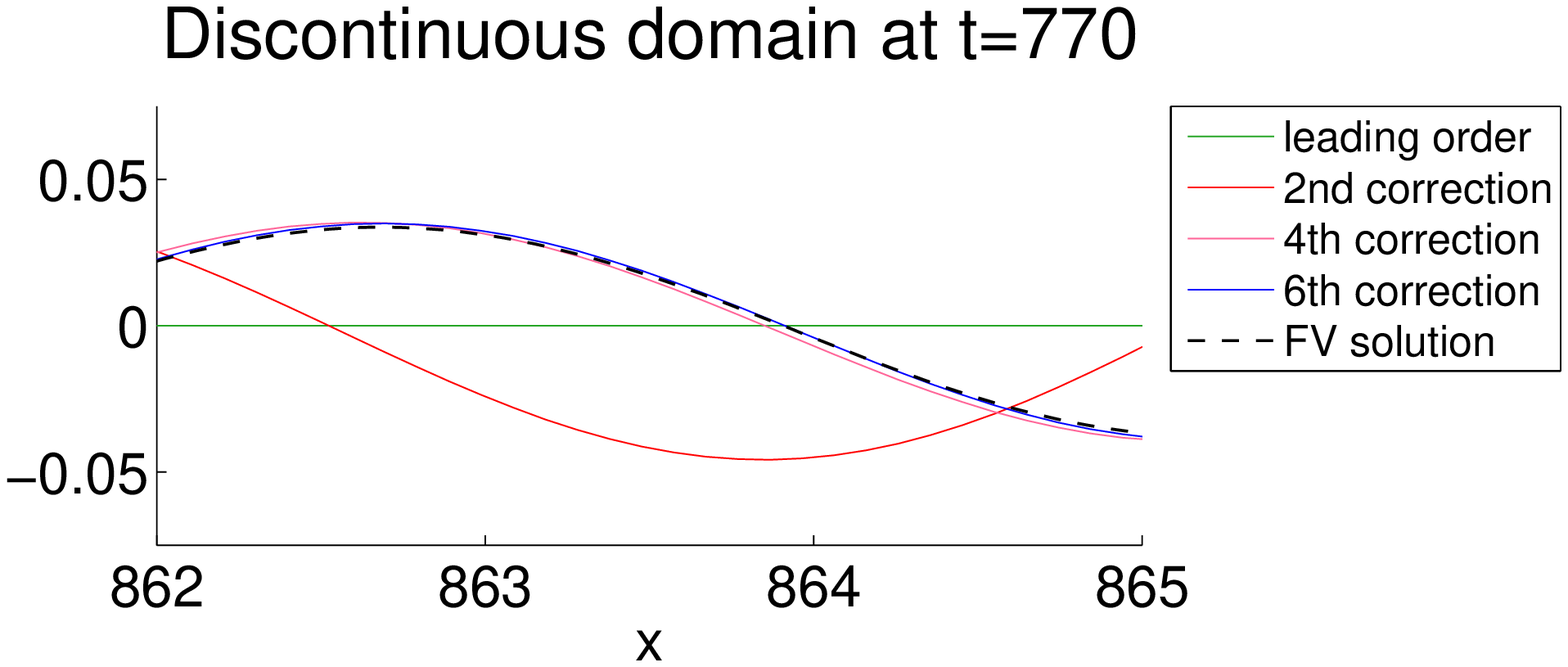}}
  
 \subfloat[Sinusoidal medium\label{fig: c-disp sin FV vs Hom}]{
  \includegraphics[scale=0.3]{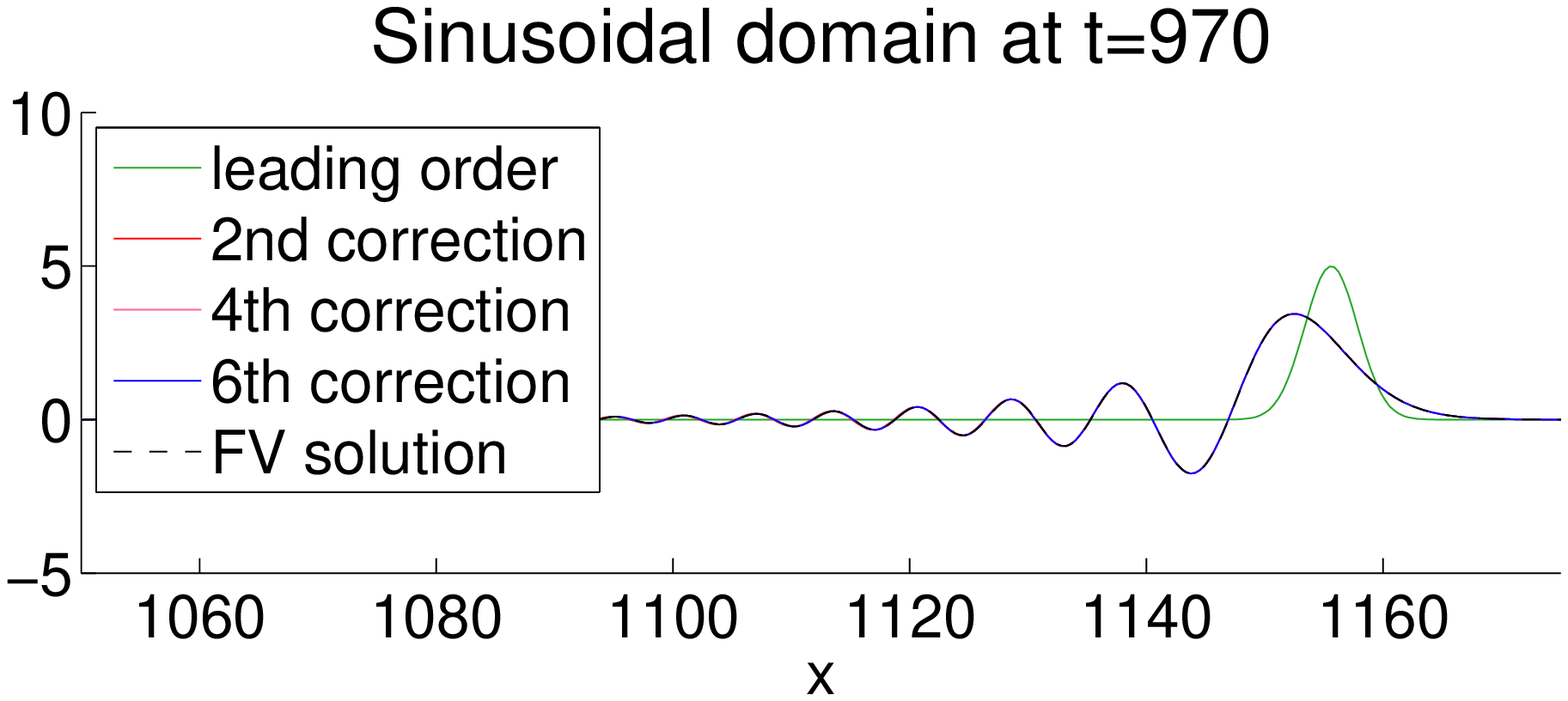}
  \includegraphics[scale=0.3]{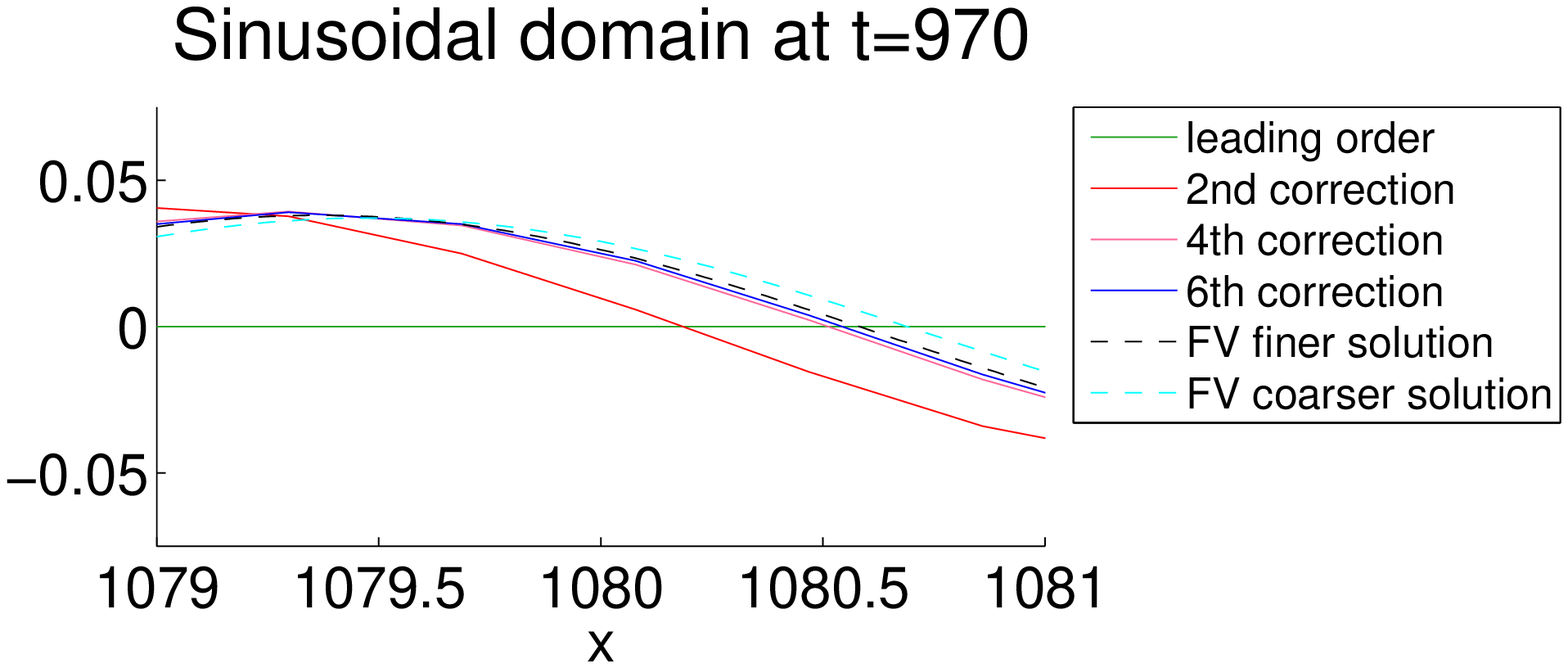}}
\par\end{centering}
\caption{Pressure from homogenized equations with different order corrections vs. $y$-averaged finite volume 
solution for (a) piecewise-constant and (b) sinusoidal media. 
On the right we show a close-up of the dispersive tail
where the differences in the homogenized corrections are more noticeable.
\label{fig: 1D wave propagation on c-dispersive media}}
\end{figure}

Figure \ref{fig: 1D wave propagation on c-dispersive media}
shows the results for the two media, after the initial pulse has traveled
a distance of more than 1000 material layers.
Homogenized solutions of differing orders are shown to demonstrate
the increasing accuracy of the high-order homogenized approximations.
These results are compared with the finite volume solution of 
the variable-coefficient wave equation \eqref{wave-equation},
which is averaged in the $y$-direction
(since the homogenized solution represents this average). 
The agreement is very good, and close examination of the dispersive tail shows
that the approximations are increasingly accurate. 

\subsection{Two-dimensional propagation}
\begin{figure}
\begin{centering}
 \subfloat[Slices of the solution along the $x$- (solid blue) and $y$-axis (solid red).
Dashed lines represent slices of the homogenized solution.\label{fig:cz-slice}]{
 \includegraphics[scale=0.3]{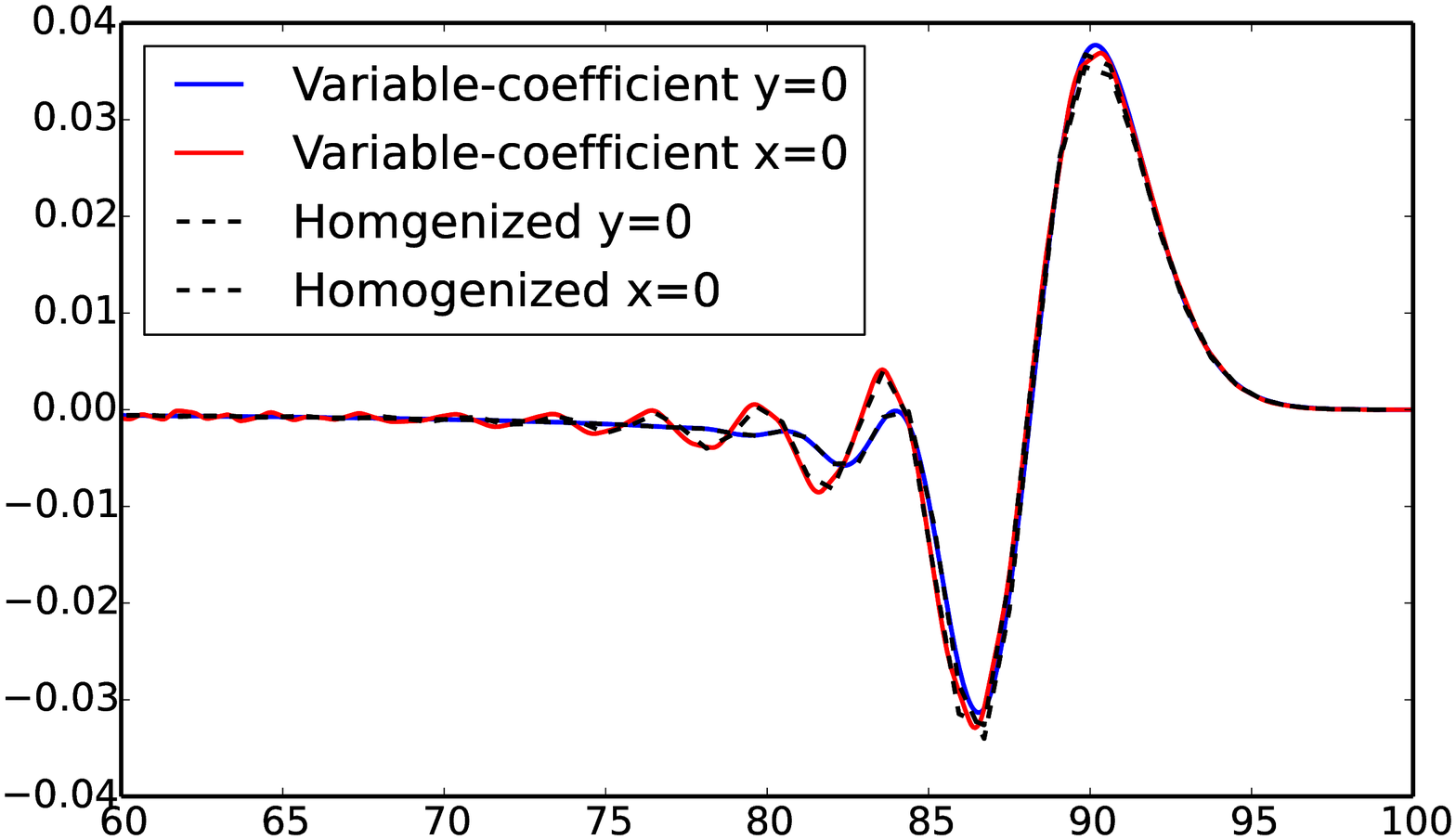}}

 \subfloat[Dispersion relation for this ``almost-isotropic'' medium.
 Speed $c=\omega(\bk)/k$ is plotted versus wavenumber.  The right plot
 shows slices along $k_x=0$ and $k_y=0$.  \label{fig:cz-disp}]{
  \includegraphics[scale=0.27]{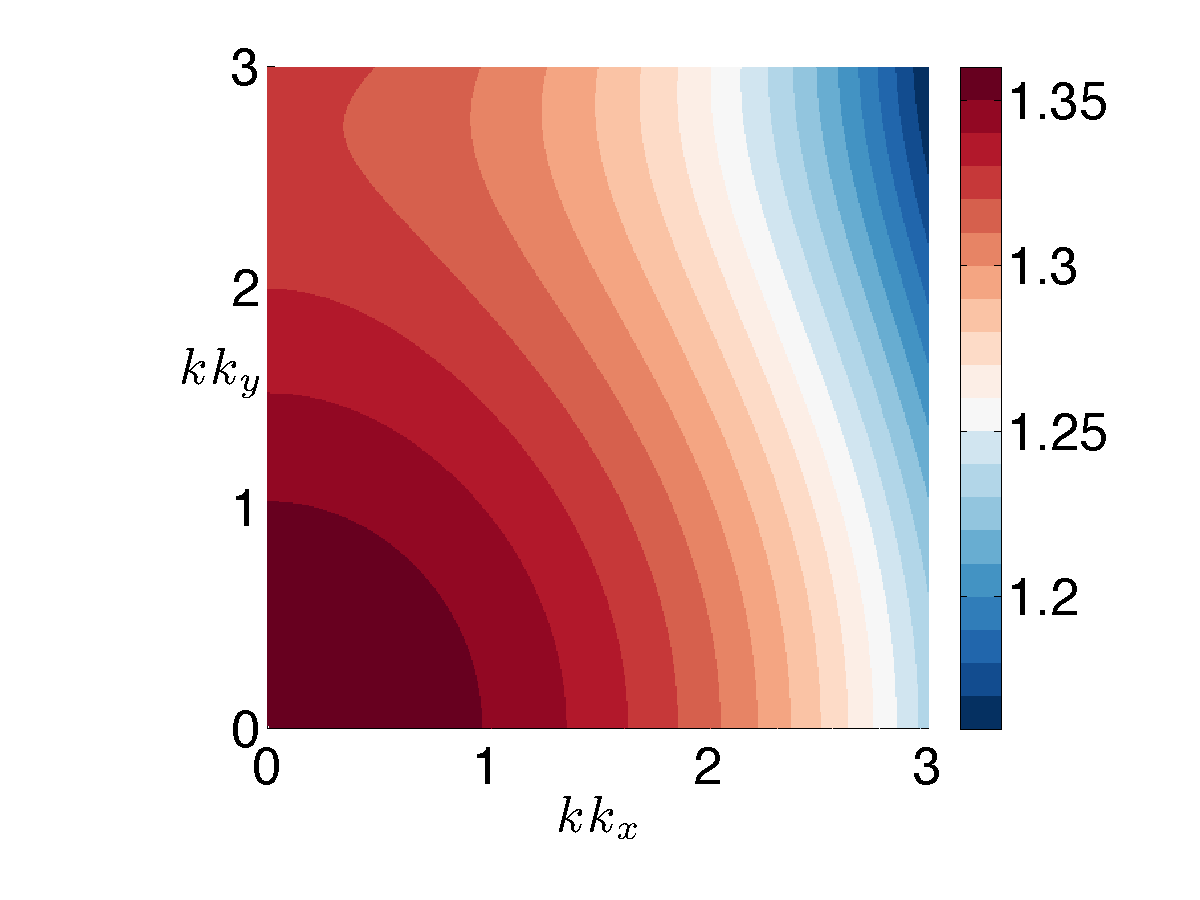}
  \includegraphics[scale=0.27]{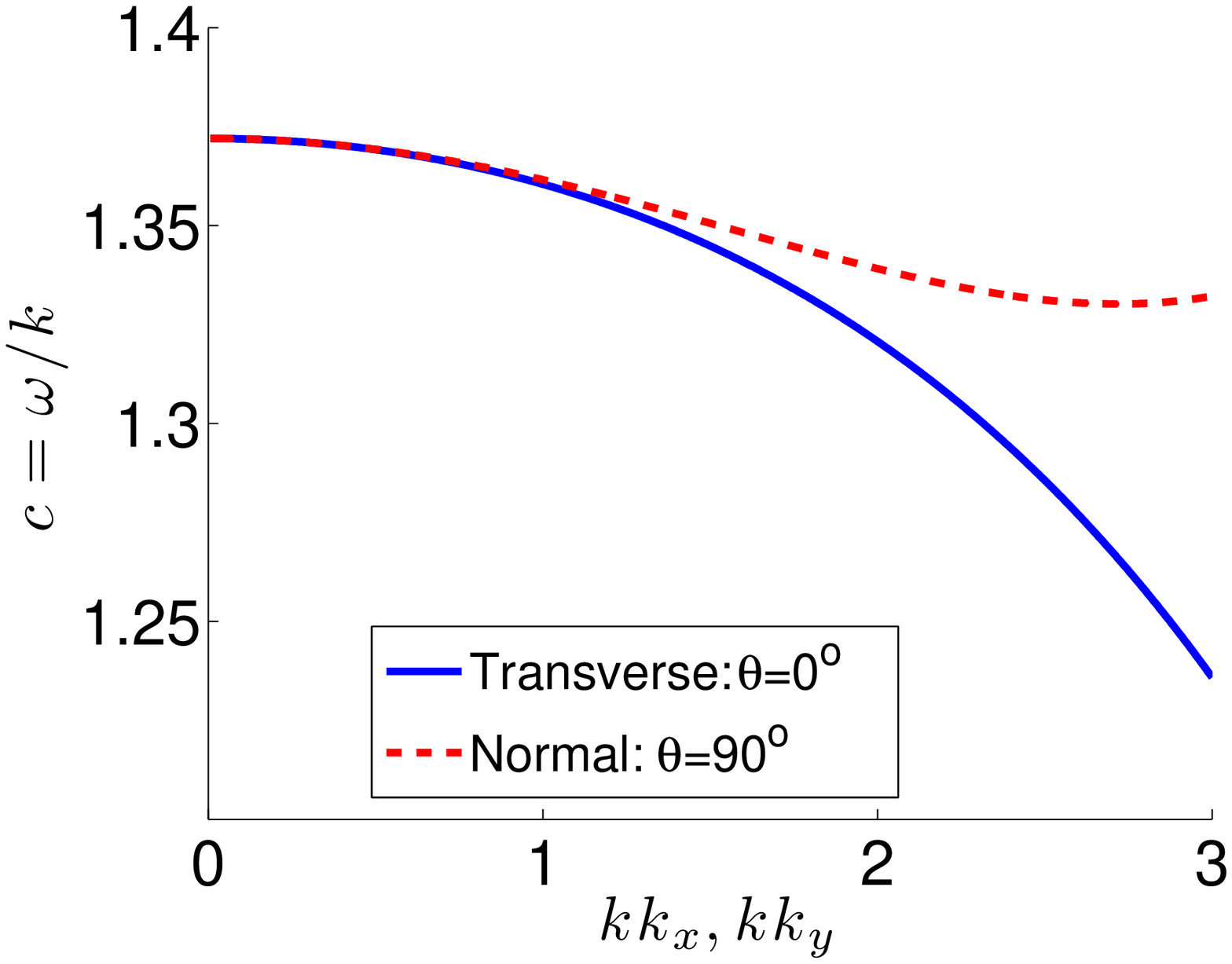}}
\par\end{centering}
\caption{2D wave propagation in a piecewise-constant medium with variable $Z$ and $c$.
We show slices along the $x$- (solid blue) and $y$-axis (solid red). We also show 
slices of the homogenized solution (dashed lines).
\label{fig:cz-wave}}
\end{figure}

We now revisit the result shown in the lower-right quadrant of 
Figure \ref{fig:quadrants_lin}.  This involves the propagation
of an initially Gaussian pressure perturbation \eqref{gaussian} in a
piecewise-constant medium \eqref{discontinuous medium} with 
\begin{align} \label{eq:almost-iso}
    K_A & = 17/2, & K_B & = 17/32, & \rho_A & = \rho_B = 1.
\end{align}
In Figure \ref{fig:cz-slice},
we show again the slices along the lines $x=0$ and $y=0$, with
solid lines indicating the variable-coefficient equation solution and dashed
lines indicating the homogenized equations solution.
Close agreement between the different approximations is evident.

Furthermore, the leading parts of the wave in the two slices are nearly
identical -- in other words, even though two very different and
direction-dependent dispersive effects are present, the medium behaves almost
isotropically.  This is the result of a special property of the chosen parameters
\eqref{eq:almost-iso}.
Due to the approximate symmetry between $Z$- and $c$-variation pointed
out in Section \ref{sec:directions}, it is possible to make the dispersion
relation isotropic up to $\order(\delta^2)$ by setting $Z_A=c_A$ and $Z_B=c_B$.
The parameters \eqref{eq:almost-iso} satisfy these conditions.
Figure \ref{fig:cz-disp} shows the dispersion relation for this medium.

\section{Discussion}
Propagation of long-wavelength waves in periodic media is known to lead to
dispersive behavior, usually resulting from small-scale reflection \cite{santosa1991}.
In this work, we have shown that macroscopic dispersion can also arise due to
microscopic diffraction, based on spatial variation of the sound speed.
We have focused on media that vary periodically in one direction and are
homogeneous in the other.
By applying homogenization techniques with carefully chosen assumptions about
which variables depend on the small scale, we have obtained a high-order accurate
effective medium approximation.  Pseudospectral solutions of the effective medium
equations agree very well with finite volume solutions of the original variable-coefficient
wave equation, even after long times.  The homogenized equations can be used to
understand the effective dispersion relation of the medium, which depends strongly
on the angle of propagation relative to the axis of variation of the medium.
By controlling the contrast in impedance and sound speed, it is possible to independently
vary the strength of dispersion in the two coordinate directions.
Interesting behaviors can be obtained by manipulating these contrasts; for instance,
we have shown that certain layered media with variable impedance and sound speed
are nevertheless effectively isotropic, up to second order.

It would be interesting to extend this study to media with periodic variation
in both $x$ and $y$.  The homogenization process becomes much more challenging
in this case, not only because there are more terms, but because it seems that
different assumptions are necessary regarding variation on the fast scale.  We
are also interested in applying our approach to nonlinear systems.  This has 
been done in \cite{quezada_diffractons} for a relatively simple system.
Extension to the shallow water equations over periodic bathymetry is the
subject of ongoing work.

The homogenization carried out here is purely formal.  It would be worthwhile to
establish precise convergence results and error bounds.  This seems to be possible
with existing techniques.

We conclude with a particularly interesting observation for which we have not
yet found a satisfactory explanation.
As pointed out in Section \ref{sec:directions}, impedance and sound speed variations
seem to play almost dual roles with respect to normal and transverse wave
propagation.  Equations \eqref{homog-z-dispersion homog system with 4 corrections} 
for normal propagation can be combined to yield the dispersive approximation
\begin{align*}
K_h^{-1} \rho_m \bar{p}_{tt} - \bar{p}_{yy} & = \delta^2(\alpha_1+\gamma_1) \bar{p}_{xxxx} + \order(\delta^4)
\end{align*}
where for the piecewise-constant medium \eqref{discontinuous medium} the
leading dispersive term coefficient is 
\begin{align*}
\alpha_1+\gamma_1 & = \frac{\left(Z_A^2 - Z_B^2\right)^2}{192 K_m^2 \rho_m^2} \lambda^2.
\end{align*}
Meanwhile, equations \eqref{homog-c-dispersion homog system with 6 corrections}
for transverse propagation can be combined to yield
\begin{align*}
K_h^{-1} \rho_h \bar{p}_{tt} - \bar{p}_{yy} & = \delta^2(\alpha_2+\beta_2) \bar{p}_{xxxx} + \order(\delta^4) 
\end{align*}
where for the piecewise-constant medium \eqref{discontinuous medium} the
leading dispersive term coefficient is 
\begin{align*}
\alpha_2+\beta_2 & = \frac{\left(c_A^2 - c_B^2\right)^2}{192 K_m^2}\rho_m\rho_h \lambda^2.
\end{align*}
The similarity is striking.  It is also striking that all dispersive coefficients
in the homogenized equations for normal propagation vanish when $Z$ is
constant, while all dispersive coefficients in the homogenized equations for transverse propagation
vanish when $c$ is constant.  We have not tried to prove that this holds in general,
but it is true for all media that we have investigated.


{\bf Acknowledgments.}  The authors thank the referees for their comments that
significantly improved this paper.  Research reported in this publication was
supported by the King Abdullah University of Science and Technology (KAUST).

\appendix
\section{Fast-variable functions} \label{sec: fast variable functions}
In this appendix we show the solution of the fast-variable functions $A$, $B$ and $C$ 
for the layered medium. They are given by:
\begin{subequations}
\begin{align}
	A(\hat{y}) &= \begin{cases} 
		\frac{\rho_h(c_A^2-c_B^2)(\lambda-4\hat{y})}{8K_m} &\text{ if } 0\leq \hat{y}\leq 1/2\lambda , \\
		\frac{\rho_h(c_B^2-c_A^2)(3\lambda-4\hat{y})}{8K_m} &\text{ if } 1/2\lambda \leq \hat{y} \leq \lambda 
			 \end{cases}, \\
	B(\hat{y}) & = \begin{cases}
		\frac{(K_A-K_B)(\lambda-4\hat{y})}{8K_m} &\text{ if } 0\leq \hat{y}\leq 1/2\lambda , \\
		-\frac{(K_A-K_B)(3\lambda-4\hat{y})}{8K_m} &\text{ if } 1/2\lambda \leq \hat{y} \leq \lambda 
				\end{cases}, \\
	C(\hat{y}) & = \begin{cases}
		-\frac{(\rho_A-\rho_B)(\lambda-4\hat{y})}{8\rho_m} &\text{ if } 0\leq \hat{y}\leq 1/2\lambda , \\
		\frac{(\rho_A-\rho_B)(3\lambda-4\hat{y})}{8\rho_m} &\text{ if } 1/2\lambda \leq \hat{y} \leq \lambda 
				\end{cases}.
\end{align}
\end{subequations}

The rest of the fast-variable functions are too cumbersome to show here. For the sinusoidal medium, 
we don't have closed form expressions; instead, we compute the coefficients of the homogenized equations numerically, 
see appendix \ref{sec: coefficients for sinusoidal medium}. 

\section{Coefficients of homogenized equations}  \label{sec: coefficients of homogenized equations}

In this appendix we give the coefficients of the homogenized systems 
\eqref{homog-combined homog system with 4 corrections}, 
\eqref{homog-c-dispersion homog system with 6 corrections} and 
\eqref{homog-z-dispersion homog system with 4 corrections}. They are:
 \begin{align*}
  \alpha_1 & = K_h\avg[ K^{-1}F], \\
  \alpha_2 & = K_h\avg[ K^{-1}H], \\
  \alpha_3 & = K_h\avg[ K^{-1}U] - K_h^2\avg[ K^{-1}F]^2, \\
  \alpha_4 & = K_h\avg[ K^{-1}W] - K_h^2\avg[ K^{-1}H]^2, \\
  \alpha_5 & = K_h\avg[ K^{-1}V] - 2K_h^2\avg[ K^{-1}F] \avg[ K^{-1}H], \\
  \alpha_6 & = K_h\avg[ K^{-1}\tilde{B}] 
  			+ K_h^3\avg[ K^{-1}H]^3
			-2K_h^2\avg[ K^{-1}H] \avg[ K^{-1}W], 
 \end{align*}
 \begin{align*}
  \beta_1 & = -\rho_h\avg[ \rho^{-1}F] & \beta_2 & = -\rho_h\avg[ \rho^{-1}H], \\
  \beta_3 & = -\rho_h\avg[ \rho^{-1}U] & \beta_4 & = -\rho_h\avg[ \rho^{-1}W], \\
  \beta_5 & = -\rho_h\avg[ \rho^{-1}V] & \beta_6 & = -\rho_h\avg[ \rho^{-1}\tilde{B}],
 \end{align*}
and 
 \begin{align*}
  \gamma_1 & = \rho_m^{-1}\avg[ \rho E], \\
  \gamma_2 & = \rho_h^{-1}\avg[ \rho D], \\
  \gamma_3 & = \rho_m^{-1}\avg[ \rho T]
  				-\rho_m^{-2}\avg[ \rho E]^2, \\
  \gamma_4 & = \rho_h^{-1}\avg[ \rho Q]
  				+\avg[ \rho^{-1} F] \avg[ \rho D]
				+ \rho_m^{-1}\avg[ \rho S]
				- \rho_m^{-1}\rho_h^{-1}\avg[ \rho E] \avg[ \rho D] \\
  \gamma_5 & = \avg[ \rho D] \avg[ \rho^{-1} H]
  				+\rho_h^{-1}\avg[ \rho R], \\
 \end{align*}
where the functions $A(\hat{y})$, $B(\hat{y})$ and $C(\hat{y})$ are defined in \eqref{A, B and C};
$D(\yhat)$, $E(\yhat)$, $F(\yhat)$ and $H(\yhat)$ are defined in \eqref{D, E, F and H};
and $Q(\yhat)$, $R(\yhat)$, $S(\yhat)$, $T(\yhat)$, $U(\yhat)$, $V(\yhat)$ and $W(\yhat)$  are given by
\begin{align*}
Q(\yhat) & = \llbracket K^{-1}K_h\left(N-CK_h\avg[ K^{-1}F ] \right) 
 			-\rho^{-1}\rho_h\left(N-C\rho_h\avg[ \rho^{-1}F ] \right) 
			-I \rrbracket,\\
R(\yhat) & = \llbracket K^{-1}K_h\left(P-CK_h\avg[ K^{-1}H ] \right) 
    -\rho^{-1}\rho_h\left(P-C\rho_h\avg[ \rho^{-1}H ] \right)-J \rrbracket,\\
S(\yhat) & = \llbracket K^{-1}K_h\left(P-CK_h\avg[ K^{-1}H ] \right)-L \rrbracket, \\ 
T(\yhat) & = \llbracket K^{-1}K_h\left(N-CK_h\avg[ K^{-1}F ] \right)-M \rrbracket, \\ 
U(\yhat) & = \llbracket \rho \rho_m^{-1}\left(M-B\rho_m^{-1}\avg[ \rho E ] \right)-N \rrbracket, \\ 
V(\yhat) & = \llbracket \rho A \avg[ \rho^{-1} F ] + \rho \rho_h^{-1}I 
    + \rho\rho_m^{-1} \left(L-B\rho_h^{-1}\avg[ \rho D ] \right)-P \rrbracket, \\ 
W(\yhat) & = \llbracket \rho A \avg[ \rho^{-1} H ]+\rho\rho_h^{-1}J \rrbracket. \\ 
\end{align*}
Finally,
\begin{align*}
\tilde{A}(\yhat) & = \llbracket K^{-1}K_h\left(W-K_h\avg[ K^{-1}W ] 
 								-K_h H \avg[ K^{-1}H ]
								+K_h^2\avg[ K^{-1}H ]^2 \right) \\
				& -\rho^{-1}\rho_h\left(W-\rho_h\avg[ \rho^{-1}W ] 
 								-\rho_h H \avg[ \rho^{-1}H ]
								+\rho_h^2\avg[ \rho^{-1}H ]^2 \right) \rrbracket,\\
\tilde{B}(\yhat) & = \llbracket \rho A \avg[ \rho^{-1} W ]
 					+ \rho J \avg[ \rho^{-1} H ]
 					+\rho\rho_h^{-1}\tilde{A} \rrbracket.
\end{align*}
In order to compute these functions, one must specify $K$ and $\rho$.
Below we give further details for the piecewise constant medium \eqref{discontinuous medium}
and the sinusoidal medium \eqref{sinusoidal medium}.
 
\subsection{Piecewise constant medium} \label{sec: coefficients for piecewise constant medium}
 
For the piecewise constant medium it is easy to obtain the 
 fast-variable functions and the coefficients in closed form; however, 
 most of them are too cumbersome to present here. Therefore, we just show 
 the coefficients of the first non-zero correction and refer to
\url{http://github.com/ketch/effective_dispersion_RR}
 for Mathematica files 
 where the rest of the coefficients can be found. The coefficients of the first non-zero 
 correction for the piecewise medium are: 
\begin{align*}
 \alpha_1 & = \frac{-\left(K_A-K_B\right)}{192K_m^2}\cdot\frac{\left(Z_A^2-Z_B^2\right)}{\rho_m} \lambda^2,\\
 \alpha_2 & = \frac{-\left(K_A-K_B\right)}{192K_m^2}\cdot\frac{\left(c_A^2-c_B^2\right)}{\rho_m^{-1}} \lambda^2,\\
 \beta_1 & = \frac{\left(\rho_A-\rho_B\right)}{192K_m}\cdot\frac{\left(Z_A^2-Z_B^2\right)}{\rho_m^2} \lambda^2,\\
 \beta_2 & = \frac{\left(\rho_A-\rho_B\right)}{192K_m}\cdot\left(c_A^2-c_B^2\right) \lambda^2,\\
 \gamma_1 & = \frac{-\left(\rho_A-\rho_B\right)}{192K_m}\cdot\frac{\left(Z_A^2-Z_B^2\right)}{\rho_m^2} \lambda^2,\\
 \gamma_2 & = \frac{-\left(\rho_A-\rho_B\right)}{192K_m}\cdot\left(c_A^2-c_B^2\right) \lambda^2.
\end{align*}

\subsection{Sinusoidal medium} \label{sec: coefficients for sinusoidal medium}
For the sinusoidal medium and for more general $y$-periodic media it is difficult  
to find closed expressions for the fast-variable functions and for the 
coefficients. Therefore, we solve the boundary value problems and compute
the coefficients numerically. Details can be found at
\url{http://github.com/ketch/effective_dispersion_RR}.
The files available there can easily be modified to produce coefficients
for other media.

The numerically computed coefficients for the sinusoidal medium are
(taking $\lambda=1$):
\begin{align*}
    \alpha_1&=2.2656\times10^{-10},  &
    \alpha_2&=-1.3208\times10^{-2},  \\
    \alpha_3&=-1.8927\times10^{-11}, &
    \alpha_4&=-1.8172\times10^{-4},  \\
    \alpha_5&=1.3398\times10^{-3}, &
    \alpha_6&=6.0711\times10^{-6}, \\
    \beta_1 &=2.9249\times10^{-4}, &
    \beta_2 &=-1.1033\times10^{-2}, \\
    \beta_3 &=-5.6345\times10^{-7},&
    \beta_4 &=-2.3474\times10^{-5}, \\
    \beta_5 &=1.1465\times10^{-3}, &
    \beta_6 &=6.9060\times10^{-6}, \\
    \gamma_1&=2.2656\times10^{-10},&
    \gamma_2&=1.2843\times10^{-2}, \\
    \gamma_3&=-1.8927\times10^{-11},&
    \gamma_4&=-1.3391\times10^{-3}, \\
    \gamma_5&=-1.6986\times10^{-4}. \\
\end{align*}

\bibliographystyle{plain}
\bibliography{linear_c-dispersion}

\end{document}